\documentclass[12pt, amscd]{amsart}
\usepackage{amscd}
\usepackage{verbatim}
\usepackage{amssymb}

\newtheorem{theorem}{Theorem}[section]
\newtheorem{proposition}[theorem]{Proposition}
\newtheorem{lemma}[theorem]{Lemma}
\newtheorem{corollary}[theorem]{Corollary}

\theoremstyle{definition}
\newtheorem{remark}[theorem]{Remark}
\newtheorem{example}[theorem]{Example}

\renewcommand{\(}{\bigl(}
\renewcommand{\)}{\bigr)}

\newcommand{\tens}{\otimes}
\newcommand{\inv}{^{-1}}
\newcommand{\xra}{\xrightarrow}

\newcommand{\iso}{\stackrel{\sim}{\to}}

\newcommand{\F}{\operatorname{F}}

\newcommand{\Pic}{\operatorname{Pic}}

\newcommand{\ch}{\operatorname{char}}

\newcommand{\cor}{\operatorname{cor}}
\newcommand{\Br}{\operatorname{Br}}
\newcommand{\Spec}{\operatorname{Spec}}

\newcommand{\Sym}{\operatorname{Sym}}

\newcommand{\Div}{\operatorname{Div}}

\newcommand{\SB}{\operatorname{SB}}

\newcommand{\gSL}{\operatorname{\mathbf{SL}}}

\newcommand{\gAut}{\operatorname{\mathbf{Aut}}}

\newcommand{\gPGL}{\operatorname{\mathbf{PGL}}}

\newcommand{\End}{\operatorname{End}}

\newcommand{\Hom}{\operatorname{Hom}}

\newcommand{\Trd}{\operatorname{Trd}}
\newcommand{\cdim}{\operatorname{cdim}}

\newcommand{\trdeg}{\operatorname{tr.deg}}

\newcommand{\G}{\mathbb{G}}
\renewcommand{\P}{\mathbb{P}}
\newcommand{\Z}{\mathbb{Z}}

\def\g{ \mathfrak{g}  }

\newcommand{\cC}{\mathcal C}

\def\F{{\overline F}}
\def\X{{\overline X}}
 \def\et{{\acute et}}

\title[Rational surfaces and canonical dimension of $\gPGL_6$] 
{Rational surfaces and canonical dimension of $\gPGL_6$}

\date{November  24, 2006}

\author
[J-L. Colliot-Th\'el\`ene] {Jean-Louis Colliot-Th\'el\`ene}

\author
[N. Karpenko] {Nikita A. Karpenko}

\author
[A. Merkurjev] {Alexander S. Merkurjev}

\address
{CNRS\\
Math\'ematiques\\
UMR  8628\\
B\^atiment 425\\
Universit\'e Paris-Sud\\
F-91405 Orsay\\
FRANCE}

\email {Jean-Louis.Colliot-Thelene {\it \`a} math.u-psud.fr}

\address
{Institut de Math\'ematiques de Jussieu\\
Universit\'e Pierre et Marie Curie - Paris 6\\
4, place Jussieu\\
F-75252 Paris CEDEX 05\\
FRANCE}

\email {karpenko {\it \`a} math.jussieu.fr}

\address
{Department of Mathematics\\
University of California\\
Los Angeles, CA 90095-1555\\
 USA}

\email {merkurev {\it at} math.ucla.edu}

\thanks{This paper is the outcome of a discussion during  a hike at Oberwolfach}

\begin{document}
\maketitle

Summary : The ``canonical dimension" of an algebraic group over a
field   by definition is the maximum of the canonical dimensions of
principal homogenous spaces under that group. Over a field of
characteristic zero, we prove that the canonical dimension of the
projective linear group $\gPGL_6$ is 3. We give two distinct proofs,
both of which rely on the birational classification of rational
surfaces over a nonclosed field. One of the proofs involves taking a
novel look at del Pezzo surfaces of degree 6.

\section{Introduction}

Let $F$ be a field and let $\cC$ be a class of field extensions of
$F$. A field $E\in\cC$ is called \emph{generic} if for any $L\in\cC$
there is an $F$-place of $E$ with values in $L$.

\begin{example}\label{am:candimvariety}
Let $X$ be a variety over $F$ and let $\cC_X$ be the class of field
extensions $L$ of $F$ such that $X(L)\neq\emptyset$. If $X$ is a
smooth irreducible variety, the field $F(X)$ is generic in $\cC$ by
\cite[Lemma 4.1]{KM}.
\end{example}

The \emph{canonical dimension} $\cdim(\cC)$ of the class $\cC$ is
the minimum of $\trdeg_F E$ over all generic fields $E\in\cC$. If
$X$ is a variety over $F$, we write $\cdim(X)$ for $\cdim(\cC_X)$
and call it the \emph{canonical dimension of $X$}. If $X$ is smooth
irreducible then by Example \ref{am:candimvariety},
\begin{equation}\label{ineq}
\cdim(X)\leq\dim X.
\end{equation}
If $X$ is smooth, proper and irreducible, the canonical dimension of
$X$ is the least dimension of a closed irreducible  subvariety $Y
\subset X$ such that there exists a rational dominant map $X
\dashrightarrow Y$ \cite[Cor .4.6]{KM}.

\begin{example}\label{severi}
Let $A$ be a central simple $F$-algebra of degree $n$. Consider the
class $\cC_A$ of all splitting fields of $A$. Let $X$ be the
Severi-Brauer variety $\SB(A)$ of right ideals in $A$ of dimension
$n$. We have $\dim X=n-1$. Since $A$ is split over a field extension
$E/F$ if and only if $X(E)\neq\emptyset$, we have $\cC_A=\cC_X$ and
therefore $\cdim(\cC_A)=\cdim(X)$.
\end{example}

Let $A$ be a central simple $F$-algebra of degree $n=q_1q_2\dots
q_r$ where the $q_i$ are powers of distinct primes. Write $A$ as a
tensor product $A_1\tens A_2\tens\dots\tens A_r$, where $A_i$ is a
central simple $F$-algebra of degree $q_i$. A field extension $E/F$
splits $A$ if and only if $E$ splits $A_i$ for all $i$. By Example
\ref{severi}, the varieties $\SB(A)$ and $Y:=\SB(A_1)\times
\SB(A_2)\times\dots\times \SB(A_r)$ have the same classes of
splitting fields and hence

\begin{equation}\label{inequality}
\cdim\SB(A)=\cdim(Y)\leq \dim(Y)=\sum_{i=1}^{r}(q_i-1)
\end{equation}
by inequality (\ref{ineq}).

It looks plausible that the inequality in (\ref{inequality}) is
actually an equality. This is proven in \cite[Th. 11.4]{BR05} in the
case when $r=1$, i.e., when $\deg(A)$ is power of a prime.

In the present paper we prove the equality in the case $n=6$.

\begin{theorem}\label{main}
Let $A$ be a division central algebra of degree $6$ over a field of
characteristic zero. Then $\cdim\SB(A)=3$.
\end{theorem}

The proof builds upon the classification of geometrically rational
surfaces. Starting from this classification, we give two independent
proofs of the theorem, each of which seems to have its own interest.
The first proof uses a novel approach to del Pezzo surfaces of
degree 6 (Section 4). The second proof involves a systematic study
of the kernel of the map from the Brauer group of a field  $F$ to
the Brauer group of the function field of a geometrically rational
surface over $F$ (Section 5).

Let $G$ be an algebraic group over $F$. The \emph{canonical
dimension of $G$} is the maximum of $\cdim(X)$ over all $G$-torsors
$X$ over all field extensions of $F$.

\begin{corollary}\label{cdpgl}
The canonical dimension of $\gPGL_6$ over a field of characteristic
zero is equal to $3$.
\end{corollary}

\begin{proof}
Isomorphism classes of $\gPGL_6$-torsors over a field extension
$E/F$ are in 1-1 correspondence with isomorphism classes of central
simple $E$-algebras of degree $6$. Moreover, if a torsor $X$
corresponds to an algebra $A$ then the classes of splitting fields
of $X$ and $A$ coincide. Therefore $\cdim(X)=\cdim\SB(A)\leq 3$.
There is a field extension $E/F$ possessing a division $E$-algebra
$A$ of degree $6$. By Theorem \ref{main}, $\cdim\SB(A)= 3$ and
therefore, $\cdim(\gPGL_6)=3$.
\end{proof}

\begin{remark}
In view of results of Berhuy and Reichstein \cite[Rem. 13.2]{BR05}
and Zainoulline \cite{Zai06},
Corollary \ref{cdpgl} completes classification of simple
groups of canonical dimension $2$ in characteristic zero. Those are
$\gSL_{3m}/{\boldsymbol{\mu}}_3$ with $m$ prime to $3$.
\end{remark}

Let $F$ be a field, $\overline F$ an algebraic closure of $F$.
An $F$-variety, or a variety over $F$, is a separated $F$-scheme of finite type.
Let $X$ be an $F$-variety. We let $\overline X =X \times_{F}{\overline F}$.

We shall use the following notation. For a variety $X$ over a field
$F$ we write $n_X$ for the \emph{index of $X$} defined as the
greatest common divisor of the degrees $[F(x):F]$ over all closed points $x\in
X$. If there exists an $F$-morphism $X \to Y$ of $F$-varieties then
$n_{Y}$ divides $n_{X}$. If $X$ is a nonempty open set of a smooth
integral quasi-projective $F$-variety $Y$ then $n_{X}=n_{Y}$ (this
may be proved by reduction to the case of a curve).  Thus if $X$ and
$Y$ are two smooth, projective, integral $F$-varieties which are
$F$-birational, then $n_{X}=n_{Y}$ (see also \cite[Rem.
6.6]{Merkurev03}).

\section{Rational curves and surfaces}

We shall need the following
\begin{theorem}\label{resolution}
Let $X$ be an integral projective variety of dimension at most $2$
over a perfect field $F$. Then there is a smooth integral projective
variety $X'$ over $F$ together with a birational morphism
$X'\to X$.
\end{theorem}

This  special case of Hironaka's theorem has been known for a long
time. In dimension 1, it is enough to normalize.  Modern proofs in
the two-dimensional case (\cite{Lip1} \cite{Lip2} \cite{Art}) handle
arbitrary excellent, noetherian two-dimensional integral schemes:
given such a scheme $X$ they produce a birational  morphism $X'
\to X$ with $X'$ regular. A regular scheme of finite
type over a perfect field $F$ is smooth over $F$.

\medskip

In this paper an integral variety $X$ over $F$ is called
\emph{rational}, resp. \emph{unirational} if there exists a
birational, resp. dominant $F$-rational map from projective space
$\P^n_{F}$ to $X$, for some integer $n$. A geometrically integral
$F$-variety $X$ is called \emph{geometrically rational}, resp.
\emph{geometrically unirational}, if there exists a birational,
resp. dominant ${\overline F}$-rational map from projective space
$\P^n_{\overline F}$ to ${\overline X}$, for some integer $n$.
Rational integral varieties are unirational. For varieties of small
dimension the converse holds under mild assumptions as the following
two well known statements show.

\begin{theorem}[L\"uroth]
\label{Lueroth}
A unirational integral curve $X$ over $F$ is
rational, i.e., $X$ is birationally isomorphic to $\P_F^1$.
\end{theorem}

\begin{theorem}[Castelnuovo]
\label{Castelnuovo}
A unirational integral surface $X$ over an algebraically closed
field field $F$ of characteristic zero is rational, i.e., $X$ is
birationally isomorphic to $\P_F^2$.
\end{theorem}

\begin{proof}
See \cite[III.2, Theorem 2.4 p.~170]{K} , or \cite{CT}. The assumption
on $\ch F$ is necessary (cf. \cite[ p.~171]{K} ).  Surfaces given by
an equation $z^p=f(x,y)$ in characteristic $p$ are unirational but
in general not rational.
\end{proof}

 \bigskip

 The following theorem has its origin in  a  paper of F.~Enriques (\cite{E}, 1897).
 The theorem as it stands was proved by V.~A.~Iskovskikh (1980) after work by
 Yu.~I.~Manin (1966, 1967).
 A proof of the theorem along the lines of modern classification theory (the cone theorem)
 was given by S.~Mori (1982).

 For a smooth $F$-variety $X$ one lets
 $K=K_{X} \in \Pic X$
 denote  the class of the canonical bundle.

A smooth   proper $F$-variety $X$ if called \emph{$F$-minimal} if
any birational $F$-morphism  from $X$ to a smooth   proper
$F$-variety is an isomorphism.  By Castelnuovo's criterion, a
smooth projective surface  over a perfect field $F$ is not
$F$-minimal if and only if $X$ contains an exceptional curve of the
first kind.

\begin{theorem}[Iskovskikh, Mori]
\label{IskoMori}

Let $X$ be a smooth, projective, geometrically integral surface over  a field $F$.
 Assume that $X$ is   geometrically rational.
 The group $\Pic X$ is free of finite type. Let $\rho$ denote its rank.
 One of the following statements holds:

(i) The surface $X$ is not $F$-minimal.

(ii) We have $\rho=2$ and $X$ is a conic bundle over a smooth conic.

(iii) We have $\rho=1$ and the anticanonical bundle $-K_{X}$ is
ample.
\end{theorem}

\begin{proof}
See \cite{I}, \cite[Thm. 2.7]{M}  and \cite[Chapter III, Section
2]{K}. See also the notes \cite{CT} (where characteristic zero is
assumed).
\end{proof}

\medskip

Smooth projective surfaces whose anticanonical bundle is ample are
known as \emph{del Pezzo surfaces}. They automatically are
geometrically rational.  Let $X/F$ be a del Pezzo surface and
$d=\deg(K_X^2)$. In particular, $n_X$ divides $d$. We have $1 \leq d
\leq 9$. For all this, see  \cite{Ma}, \cite[Chap III.3]{K},
\cite{CT}.

Over a separably closed field $F$, a del Pezzo surface is either
isomorphic to $\P^1 \times_F \P^1$ or it is obtained from $\P^2$ by
blowing up a finite set of points (at most 8, in general position).
The Picard group of $\P^2$ is $\Z h$, where $h$ is the class of a
line, and $K=-3h$. The Picard group of  $\P^1 \times_F \P^1$ is $\Z
e_{1} \oplus \Z e_{2}$, where $e_{1}$ and $e_{2}$ are the classes of
the two rulings, and $K=-2e_{1}-2e_{2}$. Given the behaviour of the
canonical class under blow-up \cite[Chap. III, Prop. 20.10]{Ma} we
therefore have:

\begin{lemma}\label{Kdivis}
Let $X/F$ be a del Pezzo surface over a separably closed field $F$.
Then one of the following mutually exclusive possibilities holds:

(i) $X$ is isomorphic to  $\P^2$.

(ii) $X$ is isomorphic to $\P^1 \times_F \P^1$.

(iii) The canonical class $K_{X}$ is not a proper mutliple of another element in $\Pic X$.
\end{lemma}

\medskip

\section{Reduction to a problem on rational surfaces}

\begin{lemma}\label{subscheme}
Let $W$ be a regular, proper, geometrically unirational variety over
a field $F$ of characteristic~$0$.  Assume that the canonical
dimension $\cdim(W)=d\leq 2$. Then there exists a geometrically
rational closed $F$-subvariety $X\subset W$ of dimension $d$ and a
dominant rational map $W \dashrightarrow X$.
\end{lemma}

\begin{proof}
By a property of canonical dimension recalled at the very beginning of this paper,
 there exist a closed irreducible
$F$-subvariety $X\subset W$ of dimension $d$ and a dominant rational
map $W\dashrightarrow X$. By assumption $W$ and therefore $X$ are
geometrically unirational. By Theorems \ref{Lueroth} and
\ref{Castelnuovo}, $X$ is a geometrically rational variety.
\end{proof}

\begin{proposition}\label{generalreduction}
Let $A$ be a division central algebra of degree 6 over  a field $F$ of characteristic zero.
Write $A=C\tens D$, where $C$ and
$D$ are central simple $F$-algebras of degree $2$ and $3$
respectively. Consider the Severi-Brauer varieties $Y=\SB(C)$ and
$Z=\SB(D)$ of dimension $1$ and $2$ respectively.
Assume $\cdim ({\rm SB(A)})\leq 2$. Then there exists a geometrically irreducible
smooth projective $F$-surface $X$ such that

$(i)$ $X$ is   $F$-minimal.

$(ii)$ $n_X$ is divisible by $6$.

$(iii)$ $X$ has a point over $F(Y\times_F Z)$.

$(iv)$  $Y \times_F Z$ has a point over $F(X)$.
\end{proposition}

\begin{proof}
Since $\cdim (Y\times_F Z)=\cdim({\rm SB(A)})   \leq 2$, then  by
Lemma \ref{subscheme}, there exist a geometrically rational closed
$F$-subvariety $X_1\subset Y\times Z$ of dimension at most $2$ and a
dominant rational map $Y\times_F Z\dashrightarrow X_1$. By Theorem
\ref{resolution}, there is a birational $F$-morphism  $X_2\to X_1$
with $X_2$ smooth and projective. Note that since $A$ is a division
algebra, we have $n_{Y\times Z}=6$. Since we have $F$-morphisms
$X_{2} \to X_{1} \to Y\times_F Z$, the numbers $n_{X_1}$ and
$n_{X_2}$ are divisible by $6$. There is a dominant rational map
$Y\times_F Z\dashrightarrow X_2$.

Suppose that $\dim X_2=1$, i.e., $X_2$ is a geometrically rational
curve. Then $X_2$ is a conic curve (twisted form of the projective
line) and $n_{X_2}$ divides $2$, a contradiction. It follows that
$X_2$ is a surface. Let $X_2 \to  X$ be a birational $F$-morphism with $X$ an
$F$-minimal smooth projective surface. Since both $X_2$ and $X$ are
smooth projective we have $n_X=n_{X_2}$.
\end{proof}

\section{Del Pezzo surfaces of degree 6}

In this Section, $F$ is an arbitrary field.

Let us first recall a few facts about del Pezzo surfaces of degree 6.
We refer to \cite{Ma} for background and proofs.

Let us first assume that $F$ is algebraically closed. A del Pezzo
surface of degree 6 is the blow-up of $\P^2$ in 3 points not on a
line.   Because $\gPGL_{3}$ acts transitively on the set of 3
noncolinear points in $\P^2$, all del Pezzo surfaces of  degree 6
are isomorphic. A concrete model is provided by the   surface $S$ in
$\P^2 \times_F \P^2$ with bihomogeneous coordinates
$[x_{0}:x_{1}:x_{2};\ y_{0}:y_{1}:y_{2}]$ defined by the system of
bihomogeneous equations $x_0y_0=x_1y_1=x_2y_2$. Projection of $S$
onto either $\P^2$ identifies $S$ with the blow-up of $\P^2$ in the
3 points $[1:0:0],\ [0:1:0]$ and $[0:0:1]$. There are 6 ``lines''
(exceptional curves of the first kind) on $S$, the inverse images
$E_{1},E_{2},E_{3}$ of the 3 points  on the first $\P^2$ and the
inverse images $F_{1},F_{2},F_{3}$ of the 3 points  on the second
$\P^2$. The configuration of these lines is that of a (regular)
hexagon : two curves $E_{i}$ do not meet, two curves $F_{i}$ do not
meet, and $(E_{i}.F_{j})=1$ if   $i\neq j$, while $(E_{i}.F_{i})=0$.

The   torus $T=(\G_{m})^3/\G_{m}$ over $F$ where $\G_{m}$ is
diagonally embedded in $(\G_{m})^3$ acts on $\P^2 \times_F \P^2$ in
the following manner:  $(t_{0},t_{1},t_{2})$ sends
$[x_{0}:x_{1}:x_{2};\ y_{0}:y_{1}:y_{2}]$ to
$[t_{0}x_{0}:t_{1}x_{1}:t_{2}x_{2};\
t_{0}^{-1}y_{0}:t_{1}^{-1}y_{1}:t_{2}^{-1}y_{2}]$. This action
induces an action on $S \subset \P^2 \times \P^2$. The torus $T$
sends each line into itself. The action of $T$ on the  complement
$U$ of the 6 lines in $S$ is faithful and transitive. If one
identifies $U$ with $T$  by the choice of a rational point in $U$,
the variety $S$ with its open set $U=T$ has the structure of a toric
variety. The symmetric group $S_{2}=\Z/2$  acts on $S \subset
\P^2\times_F \P^2$ by permuting the factors. This globally preserves
the lines, the generator of $S_{2}$ induces on the hexagon of lines
the permutation of each $E_{i}$ with each $F_{i}$, i.e. opposite
sides of the hexagon are exchanged. The group $S_{3}$ acts on $S
\subset \P^2\times_F \P^2$ by simultaneous permutation on each
factor. This globally preserves the lines. The actions of $S_{2}$
and $S_{3}$ commute. The induced action of the group $H:=S_{2}\times
S_{3}$ on the hexagon of lines realizes the automorphism group of
the hexagon.

Let the  group $H$ act on $T=(\G_{m})^3/\G_{m}$ in such a way that
the generator of $S_{2}$ sends $t \in T$ to $t^{-1}$ and $S_{3}$
acts by permutation of the factors. Let $T'$ be the semidirect
product of $T$ and $H$ with respect to this action. The above
construction yields an isomorphism from $T'$ to the algebraic group
${\gAut}(S)$ of automorphisms of the surface $S$. Indeed any
$\sigma$ in $\gAut(S)$ may be multiplied by an element of $H$ so
that the action on the hexagon becomes trivial. By general
properties of blow-ups, this implies that any of the projections $S
\to \P^2$ factorizes through the contraction $S \to \P^2$, i.e. comes
from an automorphism of $\P^2$ which respects each of the points
$[1:0:0],\ [0:1:0]$ and $[0:0:1]$. Any such automorphism is given by
an element of $T$.

Let now $F$ be an arbitrary field and $S$ a del Pezzo surface of
degree 6 over $F$. Over a separable closure ${\overline F}$ of $F$
the del Pezzo surface ${\overline S}=S \times_{F} {\overline F}$ is
split, i.e.  isomorphic to the model given above by \cite{C}.
Since the 6 lines are globally stable under the action of the Galois
group, there exists a Zariski open set $U \subset X$ whose
complement over $\overline F$ consists of the 6 lines.  The Galois
action on the 6 lines induces an automorphism of the hexagon of
lines, hence a homomorphism ${\rm Gal}{(\overline F}/F) \to
H=S_{2}\times S_{3}$. There is thus an associated \'etale quadratic
extension $K/F$ and an \'etale cubic extension $L/F$. Let $T$ be the
connected component of identity in the $F$-group $\gAut(S)$. Then
$T$ is an algebraic torus and $U$ is a principal homogeneous space
under $T$ as  they are so over $\overline F$. The group of
connected components of $\gAut(S)$ is a twisted form of $S_{2}\times
S_{3}$ : it is the $F$-group of automorphisms of the finite
$F$-scheme associated to the configuration of  the 6 lines on
${\overline F}$. The $F$-torus $T$ will be identified below (Remark
\ref{torusidentified}).

\bigskip

Let $K$ be an \'etale quadratic $F$-algebra and $B$ an Azumaya
$K$-algebra of rank $9$ over $K$ with unitary involution $\tau$
trivial on $F$ \cite[\S 2.B]{Book}. Thus $B$ is a central simple
$K$-algebra of dimension $9$ if $K$ is a field, or $B$ is isomorphic
to the product $A\times A^{op}$, where $A$ is a central simple
$F$-algebra of dimension $9$ and $A^{op}$ is the opposite algebra if
$K\simeq F\times F$.

We consider the $F$-subspace of $\tau$-symmetric elements
\[ \Sym(B,\tau)=\{b\in B\quad\text {such that}\quad \tau(b)=b\}
\]
of dimension $9$.

\begin{lemma}\label{segre}
For a right ideal $I\subset B$ of rank $3$ over $K$, the
$F$-subspace $\(I\cdot\tau I\)\cap \Sym(B,\tau)$ of $\Sym(B,\tau)$
is $1$-dimensional. The correspondence $I\mapsto \(I\cdot\tau
I\)\cap \Sym(B,\tau)$ gives rise to a closed embedding of varieties
\[
s_{B,\tau}:R_{K/F}\SB(B)\to \P\(\Sym(B,\tau)\).
\]
\end{lemma}

\begin{proof}
It is sufficient to prove the lemma in the split case, i.e., when
\[
B=\End(V)\times\End(V^*),
\]
where $V$ is a $3$-dimensional vector space over $F$, and $\tau$ is
the exchange involution $\tau(f,g^*)=(g,f^*)$ for all
$f,g\in\End(V)$ (cf. \cite[Prop. 2.14]{Book}). We identify
$\Sym(B,\tau)$ with $\End(V)$ via the embedding $f\mapsto (f,f^*)$
into $B$.

A right ideal $I$ of $B$ of rank $3$ over $K=F\times F$ is of the
form
\[
I=\Hom(V,U)\times \Hom(V^*,W^*),
\]
where $U$ and $W$ are a subspace and a factor space of $V$ of
dimension $1$ respectively. We have
\[
\tau I=\Hom(W,V)\times \Hom(U^*,V^*),
\]
and
\[
I\cdot\tau I=\Hom(W,U)\times \Hom(U^*,W^*).
\]
The $F$-space
\[
\(I\cdot\tau I\)\cap \Sym(B,\tau)=\Hom(W,U)
\]
therefore is $1$-dimensional. Under the identification of $R_{K/F}\SB(B)$ with
$\P(V)\times \P(V^*)$, the morphism $s_{B,\tau}$ takes a pair of
lines $(U,W^*)$ to $\Hom(W,U)=U\tens W^*$, i.e., it coincides with
the Segre closed embedding
\[
\P(V)\times \P(V^*)\to \P(V\tens V^*)=\P(\End(V)).\qedhere
\]
\end{proof}

Let $\Trd:B\to K$ be the reduced trace linear form. For any $x,y\in
\Sym(B,\tau)$ we have
\[
\tau\Trd(xy)=\Trd(\tau(xy))=\Trd(\tau (y)\tau
(x))=\Trd(yx)=\Trd(xy),
\]
hence $\Trd(xy)\in F$.  We therefore  have an $F$-bilinear form
${\mathfrak b}(x,y)=\Trd(xy)$ on $\Sym(B,\tau)$. The form
${\mathfrak b}$ is non-degenerate as it is so over $\overline F$.

Let $L$ be a cubic \'etale $F$-subalgebra of $B$ that is contained
in $\Sym(B,\tau)$. Write $L^\perp$ for the orthogonal complement of
$L$ in $\Sym(B,\tau)$ with respect to the form ${\mathfrak b}$. As
$L$ is \'etale, we have $L\cap L^\perp=0$. Consider the
$7$-dimensional $F$-subspace $F\oplus L^\perp$ of $\Sym(B,\tau)$ and
set
\[
S(B,\tau,L)=s_{B,\tau}\inv(\P(F\oplus L^\perp)).
\]
Thus $S(B,\tau,L)$ is a closed subvariety of $R_{K/F}\SB(B)$ (cf.
\cite{Rost96}).

An \emph{isomorphism} between two triples $(B,\tau,L)$ and
$(B',\tau',L')$ is an $F$-algebra isomorphism $f:B\iso B'$ such that
$f\circ\tau=\tau'\circ f$ and $f(L)=L'$. The automorphism group of
$(B,\tau, L)$ is a subgroup of the algebraic group $\gAut(B,\tau)$.

The construction of the scheme $S(B,\tau,L)$ being natural, an
automorphism of a triple $(B,\tau,L)$ induces an automorphism of
$S(B,\tau,L)$, i.e., we have an algebraic group homomorphism
\begin{equation}\label{isoauto}
\mu:\gAut(B,\tau,L)\to \gAut\(S(B,\tau,L)\).
\end{equation}

\begin{theorem}\label{delpezzo6}
Let $F$ be an arbitrary field. Let $B$ be a rank $9$ Azumaya algebra
with unitary involution $\tau$ over a quadratic \'etale algebra $K$
over $F$ and a cubic \'etale $F$-subalgebra $L$ of $B$ contained in
$\Sym(B,\tau)$.

(i) The variety $S(B,\tau,L)$ is a del Pezzo surface of degree $6$.

(ii)  Any del Pezzo surface of degree $6$ over $F$ is isomorphic to
$S(B,\tau,L)$ for some $B$, $\tau$ and $L$.

(iii) Two surfaces $S(B,\tau,L)$ and $S(B',\tau',L')$ are isomorphic
if and only if the triples $(B,\tau,L)$ and $(B',\tau',L')$ are
isomorphic.

(iv) The homomorphism $\mu$ is an isomorphism.

(v) The \'etale quadratic algebra  $K/F$ and the \'etale cubic
algebra $L/F$ are naturally isomorphic to the ones associated to the
Galois action on the lines of the del Pezzo surface $S(B,\tau,L)$
over ${\overline F}$.
\end{theorem}
\begin{proof}

$(i)$: We may assume that $F$ is separably closed. We claim that any
triple $(B,\tau,L)$ is isomorphic to the \emph{split triple}
$\(M_3(F)\times M_3(F),\varepsilon, F^3\)$, where:

(1) $\varepsilon(a,b)=(b^t,a^t)$, ($t$ denotes the transpose
matrix), in particular \newline $\Sym\(M_3(F)\times
M_3(F),\varepsilon\)$ consists of matrices of the shape $(a,a^t)$;

(2) $F^3$ is identified with the subalgebra of diagonal matrices in
\newline  $\Sym\(M_3(F)\times M_3(F),\varepsilon\)$, i.e. those of the
shape $(a,a)$ with $a$ diagonal;

(3) $K=F \times F \subset  M_3(F)\times M_3(F)$ is the obvious map
from $F \times F$ to the center of $M_3(F)\times M_3(F)$.

Indeed, as $K$ and $B$ are split, $(B,\tau)$ is isomorphic to
$\(M_3(F)\times M_3(F),\varepsilon\)$ by \cite[Prop. 2.14]{Book}.
Let $L'\subset M_3(F)\times M_3(F)$ be the image of the (split)
\'etale cubic subalgebra $L$ under this isomorphism. In particular,
$(B,\tau,L)\simeq \(M_3(F)\times M_3(F),\varepsilon, L'\)$. Any of
the two projections to $M_3(F)$ identifies $L'$ with a split \'etale
cubic subalgebra of $M_3(F)$. Any two split \'etale cubic
subalgebras of $M_3(F)$ are conjugate, i.e., there is an $a\in
M_3(F)^\times$ such that $aL'a\inv= F^3$. Then the conjugation by
$(a, (a\inv)^t)$ yields an isomorphism between $\(M_3(F)\times
M_3(F),\varepsilon, L'\)$ and $\(M_3(F)\times M_3(F),\varepsilon,
F^3\)$. The claim is proved.

So we may assume that $(B,\tau,L)$ is the split triple. Then
$F\oplus L^\perp$ is the space of all pairs
$(b,b^t)$
with a matrix
$b$ all diagonal elements of which are equal. Let $[x_0:x_1:x_2;\
y_0:y_1:y_2]$ be the projective coordinates in $\P^2\times_F \P^2$.
The Segre embedding $s_{B,\tau}$ takes $[x_0:x_1:x_2;\ y_0:y_1:y_2]$
to the point of $\P(M_3(F))$
given by the matrix
$(x_iy_j)_{i,j=1,2,3}$
(we here identify an element $(a,a^t) \in \Sym\(M_3(F)\times
M_3(F),\varepsilon\)$ with $a \in M_3(F)$).

Therefore $S(B,\tau,L)$ is a closed
subvariety of $\P^2\times_F\P^2$ given by the equations
$x_0y_0=x_1y_1=x_2y_2$, that is a split del Pezzo surfaces of degree
6.

\smallskip

$(iv)$: We may assume that $F$ is separably closed and hence we are
in the split situation of the proof of $(i)$. Let the torus $T$, the
semidirect product $T'$ and the split del Pezzo surface $S$ be as in
the initial discussion of del Pezzo surfaces of degree 6. We let
$T'$ act on $B$ by $F$-algebra automorphisms as follows.
The groups $T$, respectively $S_3$, act on $B$ by the formula
$x(a,b)=(xax\inv,(x\inv)^t bx^t)$, where $x$  is in $T$,
respectively is the monomial matrix corresponding to an element of
$S_3$.
The generator of
$S_2$ takes a pair $(a,b)$ to $(b,a)$. The action of $T'$ defined
this way commutes with $\tau$ and preserves $L$ elementwise and
therefore induces an algebraic group homomorphism $\varphi:T'\to
\gAut(B,\tau,L)$.
The composite map of $\varphi$ with the
homomorphism
\[
\gAut(B,\tau,L)\to\gAut(K)\times\gAut(L)=S_2\times S_3.
\]
is a surjective homomorphism which coincides with
the one described at the beginning of this section.
We claim that $\varphi$ is an isomorphism.
Let $G$
be the kernel of the above homomorphism.
It suffices to show that the restriction $\psi:T\to G$ of $\varphi$
is an isomorphism. We view $G$ as a subgroup of the connected
component $\gAut(B,\tau)^+$ of identity in $\gAut(B,\tau)$.
We have an isomorphism between $\gPGL_3$ and $\gAut(B,\tau)^+$
taking an $a$ to the conjugation by $(a,(a\inv)^t)$ (cf. \cite[\S
23]{Book}). The composite map
\[
T\xra{\psi} G\hookrightarrow\gAut(B,\tau)^+\iso \gPGL_3
\]
identifies $T$ with the maximal torus $\widetilde T$ of the classes
of diagonal matrices in $\gPGL_3$. The image of $G$ in $\gPGL_3$
coincides with the centralizer of $\widetilde T$ in $\gPGL_3$, hence
it is equal to $\widetilde T$. Thus $\psi$ is an isomorphism. The
claim is proved.

The composite map
\[
T'\xra{\varphi} \gAut(B,\tau,L)\xra{\mu} \gAut\(S(B,\tau,L)\)
\]
coincides with the isomorphism in the initial discussion of del
Pezzo surfaces of degree 6.
Therefore, $\mu$ is an isomorphism.

\smallskip

$(ii)$ and $(iii)$: By the proof of $(i)$, any triple $(B,\tau,L)$
over $\overline F$ is isomorphic to the split triple. Moreover, any
del Pezzo surface of degree 6 splits over $\overline F$. The
homomorphism $\mu$ in (\ref{isoauto}) is an isomorphism by $(iv)$,
therefore the statements follow by the standard technique in
\cite[\S 26]{Book}.

\smallskip

$(v)$: The \'etale algebras $K$ and $L$ are associated to the Galois
action on the set of 6 minimal diagonal idempotents $e_i$ and $f_i$
of the algebra $F^3\times F^3$ ($i=1,2,3$) where the $e_i$ (resp.
the $f_i$) are diagonal idempotents in $F^3\times 0$ (resp. $0\times
F^3$). The statement follows from the fact that the correspondence
$e_i\mapsto E_i$, $f_i\mapsto F_i$ establishes an isomorphism of the
$(S_2\times S_3)$-sets of minimal idempotents
$\{e_1,e_2,e_3,f_1,f_2,f_3\}$ and exceptional lines
$\{E_1,E_2,E_3,F_1,F_2,F_3\}$.
\end{proof}

\medskip

\begin{remark}\label{torusstablyrational}
With notation as in the beginning of this section, the natural exact
sequence of Galois modules
\[
0 \to {\overline F}[U]^{\times}/{\overline F}^{\times} \to \Div_{{\overline S}
\setminus {\overline U}}({\overline S}) \to\Pic{\overline S} \to \Pic{\overline U},
\]
where $\Div_{{\overline S}
\setminus {\overline U}}({\overline S})$ denotes the group of divisors of ${\overline S}$
with support on the complement of ${\overline U}$  and the first map is the divisor map,
yields the exact sequence of Galois lattices:
\[
0 \to {\hat T} \to \Z[KL/F] \to {\rm Pic} {\overline S} \to 0,
\]
which defines the 2-dimensional $F$-torus $T$ with character group $
{\hat T}= {\overline F}[U]^{\times}/{\overline F}^{\times}$. The $F$-variety $U$
is a principal homogeneous space under $T$. The 6-dimensional Galois
module $ \Z[KL/F] $ is the permutation module on the 6 lines.

Direct computation over $\overline F$ shows that there is an exact
sequence of Galois lattices
\[
0 \to {\hat T} \to \Z[KL/F] \to \Z[L/F] \oplus \Z[K/F] \to \Z \to 0.
\]
Here $ \Z[L/F]$ is the 3-dimensional permutation lattice  on the set
of opposite pairs of lines in the hexagon and $ \Z[K/F]$ is the
2-dimensional permutation lattice on the set of triangles of triples
of skew lines in the hexagon. The map  $\Z[KL/F] \to \Z[L/F]$ sends
a line to the pair it belongs to, and the map  $\Z[KL/F] \to
\Z[K/F]$ sends a line to the triangle it belongs to. The map
$\Z[L/F] \oplus \Z[K/F] \to \Z $ is the difference of the
augmentation maps. Note that this Galois homomorphism has an obvious
Galois equivariant section.

From this we conclude that there exist an isomorphism of Galois
lattices
\[
\Pic {\overline S} \oplus \Z \simeq \Z[L/F] \oplus
\Z[K/F]
\]
and  an exact sequence of $F$-tori
\[
1 \to \G_{m,F} \to R_{L/F}\G_{m} \times R_{K/F}\G_{m} \to
R_{KL/F}\G_{m} \to T \to 1.
\]
Taking $F$-points and using Hilbert's theorem 90, we conclude that
$T(F)$ is the quotient of $(KL)^{\times}$ by the subgroup spanned by $K^{\times}$
and $L^{\times}$. We also see that the $F$-torus is stably rational. More
precisely $T \times_{F} R_{K/F}\G_{m} \times_{F} R_{L/F}\G_{m}$ is
$F$-birational to $\G_{m,F}\times_{F} R_{KL/F}\G_{m}$. The $F$-torus
$T$ actually is  rational (Voskresenski\u{\i} proved that all
2-dimensional  tori are rational).
\end{remark}

\medskip

\begin{remark}\label{torusidentified}
It follows from the proof of Theorem \ref{delpezzo6} that $T$ is a
maximal $F$-torus of the connected component of the identity
$\gAut(B,\tau)^+$ of the automorphism group of the pair $(B,\tau)$.
By \cite[\S 23]{Book}, the group of $F$-points of $\gAut(B,\tau)^+$
coincides with
\[
\{b\in B^\times\ |\ b\cdot\tau(b)\in F^\times\}/K^\times.
\]
It follows that
\[
T(F)=\{x\in (KL)^\times\ |\ N_{KL/L}(x)\in F^\times\}/K^\times.
\]
We leave it to the reader to compare this description with the one produced
in the previous remark.
\end{remark}

\begin{remark}\label{ksplit}
If the quadratic algebra $K$ is split, i.e., $K=F\times F$ and
$B=A\times A^{op}$ with the switch involution $\tau$, where $A$ is a
central simple $F$-algebra of dimension $9$, the surface
$S(B,\tau,L)$ is a closed subvariety of $\SB(A)\times_F\SB(A^{op})$
and the projection $S(B,\tau,L)\to \SB(A)$ is a blow-up with center
a closed subvariety of $\SB(A)$ isomorphic to $\Spec L$. In
particular, the surface $S(B,\tau,L)$ is not minimal.
\end{remark}

\begin{lemma}\label{number}
Let $S=S(B,\tau,L)$ be a del Pezzo surfaces of degree 6. Then

   (i)
    If $n_S=6$, then $K$ and the $K$-algebra $B$ are
not split.

 (ii) If $S(F)\neq\emptyset$, then the $K$-algebra $B$ is split.
\end{lemma}

\begin{proof}
If $K$ is split then $n_S\leq 3$ by Remark \ref{ksplit}. By the same
remark, if $B$ is split then $S$ has a rational point over $K$,
hence $n_S\leq 2$. Finally, if $S$ has a rational point, then so
does $R_{K/F}\SB(B)$ as $S$ is a closed subvariety of
$R_{K/F}\SB(B)$, and therefore $B$ is split.
\end{proof}

\bigskip
We may now give our {\it first proof of Theorem} \ref{main}.
 Let notation be as in Proposition \ref{generalreduction}.

By Theorem \ref{IskoMori}, $X$ is either a conic bundle over a
smooth conic or a del Pezzo surface of degree $d=1,\dots, 9$. In the
first case, $X$ has a rational point over a field extension of
degree dividing $4$, therefore, $n_X$ divides $4$, a contradiction.
In the latter case, we have $6\ |\ n_X\ |\ d$, i.e., $d=6$ and $X$ is a
del Pezzo surface of degree 6.

By Theorem \ref{delpezzo6}, we have $X=S(B,\tau,L)$ for a rank $9$
Azumaya algebra $B$  with unitary involution $\tau$ over a quadratic
\'etale algebra $K$ over $F$ and a cubic \'etale $F$-subalgebra $L$
of $B$ contained in $\Sym(B,\tau)$. It follows from Lemma
\ref{number} (i) that $K$ and $B$ are not split. By (iii) of
Proposition \ref{generalreduction} and by Lemma \ref{number} (ii),
the $K(Y\times_F Z)$-algebra $B\tens_K K(Y\times_F  Z)$ is split.
The field extension $K(Y\times_F  Z)/K(Z)$ is the function field of
a conic over $K(Z)$ and $B\tens_K K( Z)$ is an algebra of degree $3$
over $K(Z)$, hence $B\tens_K K(Z)$ is also split. By a theorem of
Ch\^atelet (recalled below), the $K$-algebra $B$, which is not
split, is similar to $D_K$ or to $D_K^{\tens 2}$. Since $B$ carries
an involution of the second kind we have $\cor_{K/F}([B])=0$ by
\cite[Th. 3.1]{Book}. From $2[D]=\cor_{K/F}([B])$ we conclude that
$D$ and therefore $B$ is split, a contradiction.

\section{Splitting properties of geometrically rational varieties of
canonical dimension at most $2$}

In this section we study the kernel of the natural homomorphism of Brauer groups
$\Br F \to\Br F(X)$ for a geometrically unirational smooth variety $X$
of canonical dimension at most $2$.

\medskip

Let us recall the well known:

\begin{proposition}\label{Leray}
 Let $F$ be a field, $\F$ a separable closure, $\g={\rm Gal}(\F/F)$
 the absolute Galois group.
Let $X/F$ be a proper, geometrically integral variety. We then have
a natural exact sequence
$$0 \to \Pic X \to (\Pic \X)^\g \to \Br F \to \Br X,$$
where $\Br X= H^2_{\et}(X,\G_{m})$. If moreover $X/F$ is smooth,
then the map $\Br X \to \Br F(X)$ is injective, and we have the
exact sequence
$$0 \to \Pic X \to (\Pic \X)^\g \to \Br F \to \Br F(X).$$
\end{proposition}

\smallskip

We write $\Br(F(X)/F)$ for the kernel of $\Br F \to \Br F(X)$.

The following well known  result is due to F. Ch\^atelet. In
dimension $1$, i.e. for $A$ a quaternion algebra and $X$ a conic, it
goes back to Witt.

\begin{proposition}\label{chatelet}
Let $X=\SB(A)$ be the Severi-Brauer variety of~$A$. Then
$\Br(F(X)/F)$ is the subgroup of $\Br F$ generated by the class of
$A$.
\end{proposition}

\begin{proposition}\label{kerbr}
  Let $F$ be a perfect  field and $X$ a  geometrically  rational surface over $F$.
Then we have one of the following possibilities:

(i) $X$ is $F$-birational to a Severi-Brauer surface, i.e. a twisted
form of $\P^2$. Then $\Br(F(X)/F)$ is $0$ or $\Z/3$, and is spanned
by the class of a central simple algebra of degree $3$.

(ii) $X$ is $F$-birational to a twisted form of $\P^1 \times_F
\P^1$. Then $\Br(F(X)/F)$ is $0$ or $\Z/2$ (spanned by the class of
a quaternion or biquaternion algebra) or $\Z/2 \oplus \Z/2$ (spanned
by the classes of two quaternion algebras).

(iii) $X$ is $F$-birational to a conic bundle over a smooth
projective conic. Then $\Br(F(X)/F)$ is $0$ or $\Z/2$, or $\Z/2
\oplus \Z/2$, and is spanned by the classes of two quaternion
algebras.

(iv) $\Br(F(X)/F)=0$, i.e., the natural map $\Br F \to \Br F(X)$ is
injective.
\end{proposition}

\begin{proof} By resolution of singularities
(Theorem \ref{resolution}) we may assume that $X$ is smooth,
projective and $F$-minimal.

Assume $X$ is the Severi-Brauer surface $\SB(A)$  associated  to
 a central simple $F$-algebra  $A$ of index  $3$.
That $\Br(F(X)/F)$ is $0$ or $\Z/3$, and is spanned by the class of
a central simple algebra of degree $3$, follows from Proposition
\ref{chatelet}.

Assume $X$ is a twisted form of $\P^1 \times_F  \P^1$. As the
automorphism group of $\P^1 \times_F  \P^1$ is the semidirect
product of $\gPGL_2\times \gPGL_2$ with the cyclic group of order
$2$ permuting the components, we have $X=R_{K/F}(C)$ where $K/F$ is
an \'etale quadratic $F$-algebra and $C$ is a conic curve over $K$.
If $K$ is a field then by \cite[Cor. 2.12]{MT95}, we have
\[
\Br\(F(C)/F)\)=\cor_{K/F}\(\Br\(K(C)/K)\)\).
\]
By Proposition \ref{chatelet}, $\Br\(K(C)/K)\)$ is generated by the
class of a quaternion algebra over $K$ and therefore,
$\Br\(F(C)/F)\)$ is generated by the corestriction of a quaternion
algebra that is either $0$ or a quaternion algebra, or a
biquaternion algebra.

If $K=F\times F$ then $X = C \times_{F} C'$ where $C$ and $C'$ are
conics over $F$. In this case $\Br\(F(C)/F)\)$ is a quotient of
$\Z/2 \oplus \Z/2$, spanned by the classes of the quaternion
$F$-algebras associated with $C_{1}$ and $C_{2}$.

Let   $X/F$ be a conic bundle over  a conic $Y$. Then
$\Br\(F(Y)/F\)$ is $0$ or $\Z/2$, spanned by the class of the
quaternion algebra $Q$ associated to $Y$
 and the kernel of $\Br F(Y) \to \Br F(X)$ is $0$ or $\Z/2$. Thus the order of the  kernel of
 $\Br F \to \Br F(X)$ divides 4.
Let $A/F$ be a nontrivial division algebra in $\Br\(F(X)/F)\)$
different from $Q$. It suffices to show that $\Br\(F(X)/F)\)$
contains a division quaternion algebra different from $Q$. The index
of $A$ over the function field $F(Y)$ is at most 2. By the index
reduction formula \cite[Th. 1.3]{SV}, the index of one of the
$F$-algebras $A$ and $A\otimes_{F} Q$ is at most 2, i.e., one of
these two algebras is similar to a division quaternion algebra
different from $Q$.

If $X$  is not $F$-isomorphic to a twisted form of $\P^2$, to
$\P^1\times_F  \P^1$ or to a  conic bundle over a conic then
 according to Theorem \ref{IskoMori} and  Lemma \ref{Kdivis},
$X$ is  a del Pezzo surface with $\Pic X$ of rank 1 such that the
canonical class, which is in $\Pic X$, is not divisible in $\Pic
\X$. As the cokernel of the natural map $\Pic X \to (\Pic \X)^\g$ is
torsion, the group $(\Pic \X)^\g$ is free of rank $1$. Therefore
$(\Pic \X)^\g$ is generated by the canonical class and hence the map
$\Pic X \to (\Pic \X)^\g$ is an isomorphism. By Proposition
\ref{Leray}, this implies that $\Br(F(X)/F)=0$.
\end{proof}

\begin{remark}
The proof of this proposition uses Theorem \ref{IskoMori} and  Lemma
\ref{Kdivis} in a critical fashion but it requires no discussion of
del Pezzo surfaces other than Severi-Brauer surfaces and twisted
forms of $\P^1 \times_F  \P^1$. The same comment applies to Theorem
\ref{canX1or2} and Corollary \ref{3or4} hereafter, hence to the {\it
second proof of Theorem} \ref{main} given at the end of this
section.
 \end{remark}

\medskip

\begin{theorem}\label{canX1or2}
   Let $W$ be a smooth, proper, geometrically unirational variety over
a field $F$ of characteristic zero.

(i) If ${\cdim}(W)=1$ then $\Br\(F(W)/F\)$ is $0$ or $\Z/2$.

(ii) If ${\cdim}(W)=2$ then $\Br\(F(W)/F\)$ is  one of $0$,
$\Z/2$, $\Z/2 \oplus \Z/2$ or $\Z/3$. The kernel is spanned either
by    a quaternion algebra, or  two quaternion algebras, or one
biquaternion algebra, or a cubic algebra.
\end{theorem}

\begin{proof}
By Lemma \ref{subscheme}, there exists a closed geometrically
rational $F$-subvariety $X\subset W$ of dimension 1  in case (i) and
of dimension 2 in case (ii). As $\Br W$ injects into $\Br F(W)$, we
have
\[
\Br\(F(W)/F\)\subset \Br\(F(X)/F\).
\]
If the dimension of $X$ is 1, then $X$ is a smooth conic. The kernel
of  $\Br F \to \Br F(X)$ is 0 or $\Z/2$.

If the dimension of $X$ is 2, then the possibilities for
$\Br\(F(X)/F\)$ were listed in Proposition \ref{kerbr}.
\end{proof}

\smallskip

\begin{remark}  The same theorem holds if the hypothesis that $W$ is
geometrically unirational is replaced by the hypothesis that $W$ is
a geometrically rationally connected variety. These hypotheses
indeed imply that  the variety $X$ is geometrically rationally
connected. Since $X$ is of dimension at most $2$ and ${\ch}(F)=0$ this
forces $X$ to be geometrically rational.
  \end{remark}

 \medskip

 \begin{corollary}\label{3or4}
  Let $W/F$ be a smooth, proper, geometrically
unirational variety over
a field $F$ of characteristic zero.
 Assume  ${\cdim}(W) \leq 2$.  Let $A$ and $A'$ be
central division $F$-algebras. If there is an $F$-rational map from
$W$ to the product $\SB(A)\times_F  \SB(A')$ then one of the
following
  occurs:
\begin{enumerate}
  \item $A$ and $A'$ are cubic algebras.
  \item $A$ and $A'$ are quaternion or biquaternion algebras.
\end{enumerate}
\end{corollary}

\begin{proof}
If there is such a rational map, then
the classes of $A \in {\rm Ker} [\Br F \to \Br \SB(A)] $ and   $A' \in {\rm Ker}[ \Br F \to \Br \SB(A') ]$ belong to $\Br\(F(W)/F\).$
The result then follows from Theorem \ref{canX1or2}. \end{proof}

\begin{remark}  Corollary \ref{3or4} holds if the hypothesis that $W$ is
geometrically unirational is replaced by the hypothesis that $W$ is
a geometrically rationally connected variety.
 \end{remark}

 We   now give our {\it second proof of Theorem} \ref{main}.  This proof does not use Section 4.
Let notation be as in Proposition \ref{generalreduction}, so that
$Y$, resp. $Z$, is the Severi-Brauer variety attached to a
quaternion algebra, resp. to an algebra of degree 3. Assume
${\cdim}({\rm SB(A)}) \leq 2$. Then ${\cdim}(Y \times_F  Z) \leq 2$.
If we apply the above corollary \ref{3or4}  to the identity map of
$Y \times_F  Z$ then we get a contradiction. We could also combine
Proposition \ref{generalreduction} and Proposition \ref{kerbr}. It
is then clear that we here use statement (iv) of Proposition
\ref{generalreduction}, as opposed to our use of statement (iii) of
that same proposition
 in our {\it first proof} (end of Section 4)
of Theorem \ref{main}.


\end{document}